\documentclass[12pt,a4paper]{elsart}
\usepackage{mathrsfs}
\usepackage{latexsym}
\usepackage{amsmath}
\usepackage{graphicx}
\usepackage{amsfonts}
\usepackage{indentfirst}
\usepackage{amssymb}
\usepackage[english,francais]{babel}
\usepackage{epsfig,mathrsfs}
\usepackage{color}
\usepackage{cite}
\usepackage[bf,SL,BF]{subfigure}
\usepackage{extarrows}
\allowdisplaybreaks \numberwithin{equation}{section}
\newtheorem{e-proposition}[theorem]{Proposition}

\newtheorem{e-definition}[theorem]{Definition\rm}

\def\og{\leavevmode\raise.3ex\hbox{$\scriptscriptstyle\langle\!\langle$~}}
\def\fg{\leavevmode\raise.3ex\hbox{~$\!\scriptscriptstyle\,\rangle\!\rangle$}}
\date{}
\begin{document}
\selectlanguage{english}
\begin{frontmatter}
\title{Corrigendum: Recover the source and initial value simultaneously in a parabolic equation}
\noindent{2014 Inverse Problems 30 065013}
\author[author1]{Guang-Hui Zheng\corauthref{cor}},
\ead{zhgh1980@163.com}
\author[author2]{Ting Wei}
\corauth[cor]{Corresponding author.}
\address[author1]{College of Mathematics and Econometrics, Hunan University, Changsha 410082, Hunan Province, P.R. China}
\address[author2]{School of Mathematics and Statistics, Lanzhou University, Lanzhou 730000,
Key Laboratory of Applied Mathematics and Complex Systems, Gansu
Province, P. R. China}
\end{frontmatter}

In the paper above, Theorem 2 and its proof are incorrect. Because
the function $l(t)$ in (2.3) never vanishes at $t=T$, the key Lemma
\ref{lem1} is not cited correctly, so we cite Lemma 2.4 in
\cite{I(1)} as our key lemma in this corrigendum.  Due to this
consideration, we need to widen the bounded domain
$Q=\Omega\times(0,T)$ in (1.1) to $Q=\Omega\times(0,T+\delta_0)$,
where $\delta_0$ is an arbitrary fixed positive constant, i.e., we
consider the following parabolic problem
\begin{align}\label {1.1}
\begin{cases}
u_t=Au+f(x,t), &\text{ in }\ Q=\Omega\times(0,T+\delta_0),\\
\frac{\partial u}{\partial \nu_{A}}=0, &\text{ on }\ \partial\Omega\times(0,T+\delta_0),\\
u(x,0)=g(x), &\text{ in }\ \Omega
\end{cases}
\end{align}
in Theorem 2, where $A$ is a uniformly elliptic operator of second
order with $x-$dependent coefficients, and $\frac{\partial
u}{\partial \nu_{A}}$ is the conormal derivative with respect to
$u$. The admissible set is given by
\begin{align*}
U=\big\{(f,g)\big|(f,g)\in
C^{2+\gamma,\frac{2+\gamma}{2}}(\overline{Q})\times
C^{4+\gamma}(\overline{\Omega});
\|f\|_{C^{2+\gamma,\frac{2+\gamma}{2}}}(\overline{Q})+&\|g\|_{C^{4+\gamma}(\overline{\Omega})}\leq
M_0\big\},\\&\ \ \ (0<\gamma<1),
\end{align*}
and the source function $f$ in (\ref{1.1}) satisfies
\begin{align}\label{2.1}
\left|f_t(x,t)\right|\leq C_0\left|f(x,T)\right|,\ \ \ (x,t)\in
\overline{Q},
\end{align}
for some positive constant $C_0$. Then the modified function
$l(t)=t(T+\delta_0-t)$ in (\ref{2.3}) vanishes at $t=T+\delta_0$.
Moreover,
\begin{align}\label{2.3}
\rho(x,t)=\frac{e^{\lambda\psi(x)}}{l(t)},\ \ \
\theta(x,t)=\frac{e^{\lambda\psi(x)}-e^{2\lambda\|\psi(x)\|_{C({\Omega})}}}{l(t)},
\end{align}
and $\psi(x)$ is defined in \cite{I(1)} Lemma 2.3. As for the detail
proof of existence of $\psi(x)$, one can refer to Lemma 2.3 in
\cite{I(1)}. Especially, the boundary measurement is modified to
$u\big|_{\Gamma\times(T-\delta_1,T+\delta_1)}$ ($\delta_1$ will be
defined later). Compared with the common measurement
$u\big|_{\Gamma\times(0,T+\delta_0)}$, the measurement time is a
subset of whole time interval $(0,T+\delta_0)$, which is more widely
used in many applications. By doing above modification and following
\cite{I(1)}, we can actually obtain the Lipschitz stability for the
source.

All these modifications are only used in Section 2 for obtaining the
conditional stability and uniqueness, i.e. Theorem 2, Theorem 4 and
Corollary 5. Because Theorem 4 and Corollary 5 are direct results of
Theorem 2, we just focus on the corrections for Theorem 2, then the
corrections for Theorem 4 and Corollary 5 are similar. In the other
sections, the parabolic problem (\ref{1.1}) is still considered in
bounded domain $Q=\Omega\times(0,T)$ and boundary measurement is
$u\big|_{\Gamma\times(0,T)}$.

\begin{lem}\label{lem1} \cite{I(1)} There exists a number $\hat{\lambda}>0$ such that for
an arbitrary $\lambda\geq\hat{\lambda}$ we can choose $s_0(\lambda)$
such that for all $s\geq s_0(\lambda)$, the solution of parabolic
problem (\ref{1.1}) $u(f,g)\in W^{2,1}_2(Q)$ satisfies the following
inequality
\begin{align}\nonumber
&\int_Q\left((s\rho)^{p-1}\bigg(|\partial_t
u(f,g)|^2+\sum_{i,j=1}^n|\partial_i\partial_j u(f,g)|^2\right)
+(s\rho)^{p+1}|\nabla u(f,g)|^2\\
&\ \ \ \ \ \ \ \ \ \
+(s\rho)^{p+3}|u(f,g)|^2\bigg)e^{2s\theta}dxdt\leq
C\int_Q(s\rho)^{p}|f|^2e^{2s\theta}dxdt\nonumber\\
&\ \ \ \ \ \ \ \ \ \ \ \ \ \ \ \ +C\int_{\Gamma\times
(0,T+\delta_0)}\bigg((s\rho)^{p}|\partial_t
u(f,g)|^2+(s\rho)^{p+1}|\nabla u(f,g)|^2
\nonumber\\
&\ \ \ \ \ \ \ \ \ \ \ \ \ \ \ \ \ \ \ \ \ \ \ \ \ \ \ \ \ \ \ \ \ \
\ \ \ \ \ \ \ +(s\rho)^{p+3}|u(f,g)|^2\bigg)dSdt\ \ \ \
p=0,1\label{2.4}
\end{align}
where the constant $C$ depends on $\lambda$, but independent of the
large parameter $s$.
\end{lem}

\textbf{Theorem 2.} (Conditional Stability)\ For every
$\delta_1\in(0,\min\{\delta_0,T\}]$ and $(f,g)\in U$, let $u(f,g)$
be the solution of (\ref{1.1}), then we have
\begin{align*}
(1)\ \ &\|f\|_{L^2(Q)}\leq
C\|\big(u(f,g)(\cdot,T),\ u(f,g)\big)\|_{H^2(\Omega)\times {H^2(\Gamma\times(T-\delta_1,T+\delta_1))}};\\
(2)\ \ &\|g\|_{L^2(\Omega)}\leq C\left|\ln\|\big(u(f,g)(\cdot,T),\
u(f,g)\big)\|_{H^2(\Omega)\times
{H^2(\Gamma\times(T-\delta_1,T+\delta_1))}}\right|^{-1},
\end{align*}
where $C$ is a positive constant, and
\begin{align*}
\|\big(u(f,g)(\cdot,T),\ u(f,g)\big)\|&_{H^2(\Omega)\times
{H^2(\Gamma\times(T-\delta_1,T+\delta_1))}}\\
&=\left(\|u(f,g)(\cdot,T)\|_{H^2(\Omega)}^2+\|u(f,g)\|_{H^2(\Gamma\times(T-\delta_1,T+\delta_1))}^2\right)^{\frac{1}{2}}
\end{align*}

The proof of Theorem 2 is very similar to the one in \cite{I(1)},
and we correct it as follows.

\begin{pf}
(1)\ For every $\delta_1\in(0,\min\{\delta_0,T\}]$, we have $0\leq
T-\delta_1<T<T+\delta_1\leq T+\delta_0$. Then we can construct
weight functions as (\ref{2.3}) in
$Q_1=\Omega\times(T-\delta_1,T+\delta_1)$, i.e.
\begin{align*}
&l_1(t)=(t-(T-\delta_1))((T+\delta_1)-t),\\
&\rho_1(x,t)=\frac{e^{\lambda\psi(x)}}{l_1(t)},\\
&\theta_1(x,t)=\frac{e^{\lambda\psi(x)}-e^{2\lambda\|\psi(x)\|_{C(\overline{\Omega})}}}
{l_1(t)},
\end{align*}
for $(x,t)\in Q_1$, where $\psi(x)$ is defined in (\ref{2.3}).
Similar to \cite{I(1)}, by the time transform
$\tilde{t}=t-T+\delta_1$, we can change $Q_1$ into
$\tilde{Q}_1=\{(x,\tilde{t})|(x,\tilde{t})\in\Omega\times(0,2\delta_1)\}$,
and change $Q=\{(x,t)|(x,t)\in\Omega\times(0,T+\delta_0)\}$ into
$\tilde{Q}=\{(x,\tilde{t})|(x,\tilde{t})\in\Omega\times(-T+\delta_1,
\delta_0+\delta_1)\}$. We focus on the domain $\tilde{Q}_1$, using
the transform above, the weight functions become into
\begin{align*}
&\tilde{l}_1(\tilde{t})=\tilde{t}(2\delta_1-\tilde{t}),\\
&\tilde{\rho}_1(x,\tilde{t})=\frac{e^{\lambda\psi(x)}}{\tilde{l}_1(\tilde{t})},\\
&\tilde{\theta}_1(x,\tilde{t})=\frac{e^{\lambda\psi(x)}-e^{2\lambda\|\psi(x)\|_{C(\overline{\Omega})}}}
{\tilde{l}_1(\tilde{t})},
\end{align*}
for $(x,\tilde{t})\in \tilde{Q}_1$. Setting
$\tilde{u}(x,\tilde{t}):=u(x,\tilde{t}+T-\delta_1)$ and
$\tilde{f}(x,\tilde{t}):=f(x,\tilde{t}+T-\delta_1)$, then
$\tilde{u}$ in $\tilde{Q}_1$ satisfies the following parabolic
equation
\begin{align}\label{}
\begin{cases}
\tilde{u}_{\tilde{t}}=A\tilde{u}+\tilde{f}(x,\tilde{t}), &\text{ in }\ \tilde{Q}_1,\\
\frac{\partial \tilde{u}}{\partial \nu_{A}}=0, &\text{ on }\ \partial\Omega\times(0,2\delta_1),\\
\tilde{u}(x,0)=u(x,T-\delta_1), &\text{ in }\ \Omega.
\end{cases}
\end{align}
Since $(f,g)\in U$, the solution of parabolic problem (\ref{1.1})
$u\in C^{4+\gamma,\frac{4+\gamma}{2}}(\overline{Q})$. We define
$v:=\tilde{u}_{\tilde{t}}$, then $v_{\tilde{t}}$, $Av$ exist and $v$
satisfies
\begin{align}\label{2.5}
\begin{cases}
v_{\tilde{t}}=Av+\tilde{f}_{\tilde{t}}(x,\tilde{t}), &\text{ in }\ \tilde{Q}_1,\\
\frac{\partial v}{\partial \nu_{A}}=0, &\text{ on }\ \partial\Omega\times(0,2\delta_1),\\
v(x,0)=\tilde{u}_{\tilde{t}}(x,0), &\text{ in }\ \Omega.
\end{cases}
\end{align}
Owing to $\tilde{f}_{\tilde{t}}\in
C^{\gamma,\frac{\gamma}{2}}(\overline{\tilde{Q}_1})$ and
$\tilde{u}_{\tilde{t}}(x,0)\in C^{2+\gamma}(\overline{\Omega})$, we
see the solution of (\ref{2.5}) $v\in
C^{2+\gamma,\frac{2+\gamma}{2}}(\overline{\tilde{Q}_1})$. Noting
that the weight functions in $\tilde{Q}_1$ are consistent with the
ones in $Q$, we can still use the Carleman estimate (\ref{2.4}) in
$\tilde{Q}_1$. Thus, from Lemma \ref{lem1} with $p=0$, we get the
Carleman inequality for $v$, that is, there exists $\hat{\lambda}>0$
such that for $\lambda=\hat{\lambda}$ we can choose
$s_0(\hat{\lambda})$ such that for all $s\geq s_0(\hat{\lambda})$,
$v$ satisfies
\begin{align}\nonumber
&\int_{\tilde{Q}_1}\left(\frac{1}{s\tilde{\rho}_1}\left(|\partial_{\tilde{t}}
v|^2+\sum_{i,j=1}^n|\partial_i\partial_j v|^2\right)
+s\tilde{\rho}_1|\nabla
v|^2+s^3\tilde{\rho}_1^3|v|^2\right)e^{2s\tilde{\theta}_1}dxd\tilde{t}\\\leq
&C\int_{\tilde{Q}_1}|\tilde{f}_{\tilde{t}}(x,\tilde{t})|
^2e^{2s\tilde{\theta}_1}dxd\tilde{t}+C\int_{\Gamma\times(0,2\delta_1)}\left(|\partial_{\tilde{t}}
v|^2+s\tilde{\rho}_1|\nabla
v|^2+s^3\tilde{\rho}_1^3|v|^2\right)dSd\tilde{t},\label{2.6}
\end{align}
where we set $\lambda=\hat{\lambda}$ in $\tilde{\rho}_1$,
$\tilde{\theta}_1$, and throughout this section, $C$ always denotes
a positive generic constant which depends on $\hat{\lambda}$, but
independent of large parameter $s$.

In particular, from above time linear transform, we find the
measured time $t=T$ is changed into $\tilde{t}=\delta_1$. Therefore,
in $\tilde{Q}$, the condition (\ref{2.1}) becomes into
\begin{align}\label{3.4}
\left|\tilde{f}_{\tilde{t}}(x,\tilde{t})\right|\leq
C_0\left|\tilde{f}(x,\delta_1)\right|,\ \ \ (x,\tilde{t})\in
\overline{\tilde{Q}}.
\end{align}
Since
$\tilde{\theta}_1(x,\tilde{t})\leq\tilde{\theta}_1(x,\delta_1)$, for
$(x,\tilde{t})\in \tilde{Q}_1$, and from the condition (\ref{3.4}),
then (\ref{2.6}) yields
\begin{align}\nonumber
&\int_{\tilde{Q}_1}\left(\frac{1}{s\tilde{\rho}_1}\left(|\partial_{\tilde{t}}
v|^2+\sum_{i,j=1}^n|\partial_i\partial_j v|^2\right)
+s\tilde{\rho}_1|\nabla
v|^2+s^3\tilde{\rho}_1^3|v|^2\right)e^{2s\tilde{\theta}_1}dxd\tilde{t}\\\leq
&C\int_\Omega\left|\tilde{f}(x,\delta_1)\right|^2
e^{2s\tilde{\theta}_1(x,\delta_1)}dx+C\int_{\Gamma\times(0,2\delta_1)}\left(|\partial_{\tilde{t}}
v|^2+s\tilde{\rho}_1|\nabla
v|^2+s^3\tilde{\rho}_1^3|v|^2\right)dSd\tilde{t},\nonumber\\
&\forall s>s_0(\hat{\lambda}).\label{2.7}
\end{align}
According to
$v(x,\delta_1)=A_{\delta_1}\tilde{u}+\tilde{f}(x,\delta_1)$, where
\begin{align}\label{at}
A_{\delta_1}q=\sum_{i,j=1}^n\frac{\partial}{\partial
x_i}\left(a_{ij}(x)\frac{\partial q}{\partial
x_j}(x,\delta_1)\right)+\sum_{i=1}^n b_{i}(x)\frac{\partial
q}{\partial x_i}(x,\delta_1)+c(x)q(x,\delta_1),
\end{align}
we have
\begin{align}\label{2.8}
\int_\Omega
s\left|\tilde{f}(x,\delta_1)\right|^2e^{2s\tilde{\theta}_1(x,\delta_1)}dx\leq2\int_\Omega
s\left|v(x,\delta_1)\right|^2&e^{2s\tilde{\theta}_1(x,\delta_1)}dx\nonumber\\&+2\int_\Omega
s\left|A_{\delta_1}\tilde{u}\right|^2e^{2s\tilde{\theta}_1(x,\delta_1)}dx.
\end{align}
Thanks to the construction of $\tilde{\theta}_1(x,\tilde{t})$, the
following inequality is hold
\begin{align}
\tilde{\theta}_1(x,\tilde{t})\leq-\frac{\left(e^{2\hat{\lambda}\|\psi(x)\|_{C(\overline{\Omega})}}-e^{\hat{\lambda}\|\psi(x)\|
_{C(\overline{\Omega})}}\right)}{\delta_1^2},\ \ \ (x,\tilde{t})\in
\overline{\tilde{Q}_1}.
\end{align}
writing
$M:=\frac{\left(e^{2\hat{\lambda}\|\psi(x)\|_{C(\overline{\Omega})}}-e^{\hat{\lambda}\|\psi(x)\|
_{C(\overline{\Omega})}}\right)}{\delta_1^2}>0$, by utilizing
$v^2(x,\tilde{t})e^{2s\tilde{\theta}_1(x,\tilde{t})}\rightarrow0\
(\tilde{t}\rightarrow0^+)$ and
$|\partial_{\tilde{t}}\tilde{\theta}_1|\leq C \tilde{\rho}_1^2$ in
$\overline{\tilde{Q}_1}$, (\ref{2.8}) implies
\begin{align*}
\int_\Omega
s&\left|\tilde{f}(x,\delta_1)\right|^2e^{2s\tilde{\theta}_1(x,\delta_1)}dx\\
&\leq2\int_0^{\delta_1}\frac{\partial}{\partial
\tilde{t}}\left(\int_\Omega s
v^2(x,\tilde{t})e^{2s\tilde{\theta}_1(x,\tilde{t})}dx\right)d\tilde{t}+Cse^{-2sM}\left\|\tilde{u}(x,\delta_1)\right\|_{H^2(\Omega)}^2\\
&=2\int_{\Omega\times(0,\delta_1)}\left(2sv(\partial_{\tilde{t}}
v)+2s^2(\partial_{\tilde{t}}\tilde{\theta}_1)v^2\right)e^{2s\tilde{\theta}_1(x,\tilde{t})}dxd\tilde{t}
+Cse^{-2sM}\left\|\tilde{u}(x,\delta_1)\right\|_{H^2(\Omega)}^2\\
&\leq2\int_{\Omega\times(0,\delta_1)}\left(\frac{2}{\sqrt{s\tilde{\rho}_1}}(\partial_{\tilde{t}}
v)e^{s\tilde{\theta}_1}\right)\left(s\sqrt{s\tilde{\rho}_1}v
e^{s\tilde{\theta}_1}\right)dxd\tilde{t}+C\int_{\Omega\times(0,\delta_1)}s^2\tilde{\rho}_1^2v^2e^{2s\tilde{\theta}_1}dxd\tilde{t}\\
&\ \ \ \ \ \ \ +Cse^{-2sM}\left\|\tilde{u}(x,\delta_1)\right\|_{H^2(\Omega)}^2\\
&\leq2\int_{\Omega\times(0,\delta_1)}\frac{1}{s\tilde{\rho}_1}|\partial_{\tilde{t}}
v|^2e^{2s\tilde{\theta}_1}dxd\tilde{t}+2\int_{\Omega\times(0,\delta_1)}s^3\tilde{\rho}_1
v^2e^{2s\tilde{\theta}_1}dxd\tilde{t}\\&\ \ \ \ \ \ \ \ \ \ \ \ \ \
\ \ \ \ \ \
+C\int_{\Omega\times(0,\delta_1)}s^2\tilde{\rho}_1^2v^2e^{2s\tilde{\theta}_1}dxd\tilde{t}
+Cse^{-2sM}\left\|\tilde{u}(x,\delta_1)\right\|_{H^2(\Omega)}^2.
\end{align*}
Compare the last inequality with (\ref{2.7}), it follows that
\begin{align}\nonumber
(s-C)\int_\Omega
\left|\tilde{f}(x,\delta_1)\right|^2e^{2s\tilde{\theta}_1(x,\delta_1)}&dx
\leq C\int_{\Gamma\times(0,2\delta_1)}\left(|\partial_{\tilde{t}}
v|^2+s\tilde{\rho}_1|\nabla
v|^2+s^3\tilde{\rho}_1^3|v|^2\right)dSd\tilde{t}\\&+Cse^{-2sM}\left\|\tilde{u}(x,\delta_1)\right\|_{H^2(\Omega)}^2,\
\ \ \forall s>s_0(\hat{\lambda}).\label{2.10}
\end{align}
On the other hand, in term of (\ref{3.4}) in $\tilde{Q}$, we find
\begin{align}\nonumber
\int_{\tilde{Q}}|\tilde{f}(x,\tilde{t})|^2e^{2s\tilde{\theta}_1(x,\delta_1)}dxd\tilde{t}=&\int_{\tilde{Q}}\left|-\int_{\tilde{t}}^{\delta_1}
\tilde{f}_\xi(x,\xi)d\xi+\tilde{f}(x,\delta_1)\right|^2e^{2s\tilde{\theta}_1(x,\delta_1)}dxd\tilde{t}\nonumber
\\\leq&\int_{\tilde{Q}}\left(\left|\int_{\tilde{t}}^{\delta_1} |\tilde{f}_\xi(x,\xi)|d\xi\right|+\left|\tilde{f}(x,\delta_1)\right|\right)^2
e^{2s\tilde{\theta}_1(x,\delta_1)}dxd\tilde{t}\nonumber
\\\leq& C\int_\Omega\left|\tilde{f}(x,\delta_1)\right|^2
e^{2s\tilde{\theta}_1(x,\delta_1)}dx.\label{2.11}
\end{align}
From (\ref{2.10}) and (\ref{2.11}), it follows that
\begin{align}\nonumber
(s-C)\int_{\tilde{Q}}
|\tilde{f}(x,\tilde{t})|^2e^{2s\tilde{\theta}_1(x,\delta_1)}&dxd\tilde{t}
\leq C\int_{\Gamma\times(0,2\delta_1)}\left(|\partial_{\tilde{t}}
v|^2+s\tilde{\rho}_1|\nabla
v|^2+s^3\tilde{\rho}_1^3|v|^2\right)dSd\tilde{t}\\&+Cse^{-2sM}\left\|\tilde{u}(x,\delta_1)\right\|_{H^2(\Omega)}^2,\
\ \ \forall s>s_0(\hat{\lambda}).\label{2.12}
\end{align}
Furthermore, setting
$s_1(\hat{\lambda}):=\max\{s_0(\hat{\lambda}),2C\}$, and we obtain
\begin{align}\nonumber
\frac{1}{2}\int_{\tilde{Q}}
|\tilde{f}(x,\tilde{t})|^2e^{2s\tilde{\theta}_1(x,\delta_1)}dxd\tilde{t}
&\leq
\frac{1}{2}\int_{\Gamma\times(0,2\delta_1)}\left(|\partial_{\tilde{t}}
v|^2+s\tilde{\rho}_1|\nabla
v|^2+s^3\tilde{\rho}_1^3|v|^2\right)dSd\tilde{t}\\&\ \ \ \ \ \
+Ce^{-2sM}\left\|\tilde{u}(x,\delta_1)\right\|_{H^2(\Omega)}^2,\ \ \
\forall s>s_1(\hat{\lambda}).\label{}
\end{align}
In view of the continuity of $\tilde{\theta}_1(x,\delta_1)$, we see
there exist a positive constant $c_1(\hat{\lambda})$ such that
$\tilde{\theta}_1(x,\delta_1)\geq-c_1(\hat{\lambda})$, $\forall
x\in\overline{\Omega}$, and so
\begin{align}\nonumber
\int_{\tilde{Q}} |\tilde{f}(x,\tilde{t})|^2dxd\tilde{t} &\leq
Cs^3e^{2c_1(\hat{\lambda})s}\int_{\Gamma\times(0,2\delta_1)}\left(|\partial_{\tilde{t}}
v|^2+|\nabla v|^2+|v|^2\right)dSd\tilde{t}\\&\ \ \ \ \ \
+Ce^{2(c_1(\hat{\lambda})-M)s}\left\|\tilde{u}(x,\delta_1)\right\|_{H^2(\Omega)}^2,\
\ \ \forall s>s_1(\hat{\lambda}).\label{3.12}
\end{align}
Next we fix $s$ in the right-hand side of (\ref{3.12}), it concludes
\begin{align}
&\left(\int_{\tilde{Q}}
|\tilde{f}(x,\tilde{t})|^2dxd\tilde{t}\right)^{\frac{1}{2}}\leq
C\left\|\big(\tilde{u}(\cdot,\delta_1),\
\tilde{u}(\cdot,\cdot)\big)\right\|_{H^2(\Omega)\times
{H^2(\Gamma\times(0,2\delta_1))}}.\label{2.15}
\end{align}
Hence, noting the time inverse transform, we convert back to the
$t$-variable and obtain (1).

(2)\ We directly write $\vartheta:=u_t$ in (\ref{1.1}) and have
\begin{align}\label{2.17}
\begin{cases}
\vartheta_t=A\vartheta+f_t(x,t), &\text{ in }\ Q,\\
\frac{\partial \vartheta}{\partial \nu_{A}}=0, &\text{ on }\ \partial\Omega\times(0,T+\delta_0),\\
\vartheta(x,T)=A_{T} u+f(x,T), &\text{ in }\ \Omega,
\end{cases}
\end{align}
where the operator $A_T$ is defined as (\ref{at}). We decompose
(\ref{2.17}) as follows,
\begin{align}
\begin{cases}
w_t=Aw+f_t(x,t), &\text{ in }\ Q,\\
\frac{\partial w}{\partial \nu_{A}}=0, &\text{ on }\ \partial\Omega\times(0,T+\delta_0),\\
w(x,0)=0, &\text{ in }\ \Omega,
\end{cases}
\end{align}
and
\begin{align}\label{2.19}
\begin{cases}
z_t=Az, &\text{ in }\ Q,\\
\frac{\partial z}{\partial \nu_{A}}=0, &\text{ on }\ \partial\Omega\times(0,T+\delta_0),\\
z(x,T)=A_{T} u+f(x,T)-w(x,T), &\text{ in }\ \Omega.
\end{cases}
\end{align}
Clearly, $\vartheta=w+z$, and then $\vartheta(x,0)=z(x,0)$, for all
$x\in\Omega$. Similar to the solution of (\ref{2.5}), we find
$\vartheta\in C^{2+\gamma,\frac{2+\gamma}{2}}(\overline{Q})$.
Consequently,
\begin{align*}
\|z(\cdot,0)\|_{L^{\infty}(\Omega)}=\|\vartheta(\cdot,0)\|_{L^{\infty}(\Omega)}
\leq C\|\vartheta\|_{C^{2+\gamma,
\frac{2+\gamma}{2}}(\overline{Q})}\leq CM_0.
\end{align*}
Applying the well-know result (For example \cite{I(4)}), we have
\begin{align*}
\|z(\cdot,t)\|_{L^{2}(\Omega)}\leq
(CM_0)^{1-\frac{t}{T}}\cdot\|z(\cdot,T)\|_{L^{2}(\Omega)}^{\frac{t}{T}},\
\ \ \ \ t\in[0,T].
\end{align*}
Furthermore, by the semigroup theory (See \cite{P(5)}), we get
\begin{align*}
w(\cdot,t)=w(t)=\int_{0}^tS(t-\tau)f_\tau(\tau)d\tau,
\end{align*}
where $S(t)$, $t\geq0$ is the $C_0$-semigroup generated by $A$, and
\begin{align*}
D(A)=\left\{u\in L^2(\Omega)\mid Au\in L^2(\Omega),\frac{\partial
u}{\partial \nu_{A}}\big| _{\partial\Omega}=0\right\}.
\end{align*}
By the property of $C_0$-semigroup and condition (\ref{2.1}), it
follows that
\begin{align*}
\|w(\cdot,t)\|_{L^2(\Omega)}=&\int_{0}^t\|S(t-\tau)\|\cdot\|f_\tau(\cdot,\tau)\|_{L^2(\Omega)}d\tau\\
\leq&\int_{0}^tC\|f(\cdot,T)\|_{L^2(\Omega)}d\tau\\
\leq& C\|f(\cdot,T)\|_{L^2(\Omega)},
\end{align*}
for all $t\in[0,T]$. Employing (\ref{2.19}), we can estimate
\begin{align*}
\|\vartheta(\cdot,t)\|_{L^2(\Omega)}\leq&\|z(\cdot,t)\|_{L^2(\Omega)}+\|w(\cdot,t)\|_{L^2(\Omega)}\\
\leq&C\|z(\cdot,T)\|_{L^2(\Omega)}^{\frac{t}{T}}+C\|f(\cdot,T)\|_{L^2(\Omega)}\\
\leq&C(\|u(\cdot,T)\|_{H^2(\Omega)}+\|f(\cdot,T)\|_{L^2(\Omega)})^{\frac{t}{T}}
+C\|f(\cdot,T)\|_{L^2(\Omega)}.
\end{align*}
Therefore, utilize the $L^2-$ estimation of $u(\cdot,T)$ in (1)
(such as (\ref{2.10})), and note that $\|u(\cdot,T)\|_{H^2(\Omega)}$
and $\|u\|_{H^2(\Gamma\times(T-\delta_1,T+\delta_1))}$ will be small
enough, it implies
\begin{align*}
\|g\|_{L^2(\Omega)}=&\|u(\cdot,0)\|_{L^2(\Omega)}\\
=&\|-\int_{0}^{T} \vartheta(\cdot,\tau)d\tau+u(\cdot,T)\|_{L^2(\Omega)}\\
\leq&C\int_{0}^{T}(\|u(\cdot,T)\|_{H^2(\Omega)}+\|f(\cdot,T)\|_{L^2(\Omega)})
^{\frac{\tau}{T}}d\tau\\
&+C\|f(\cdot,T)\|_{L^2(\Omega)}+\|u(\cdot,T)\|_{L^2(\Omega)}\\
\leq&C\frac{\left|1-\|\big(u(f,g)(\cdot,T),\
u(f,g)\big)\|_{H^2(\Omega)\times
{H^2(\Gamma\times(T-\delta_1,T+\delta_1))}}\right|}{\left|\ln
\left(\|\big(u(f,g)(\cdot,T),\ u(f,g)\big)\|_{H^2(\Omega)\times {H^2(\Gamma\times(T-\delta_1,T+\delta_1))}}\right)\right|}\\
&+C\|\big(u(f,g)(\cdot,T),\ u(f,g)\big)\|_{H^2(\Omega)\times {H^2(\Gamma\times(T-\delta_1,T+\delta_1))}}\\
\leq&C\left|\ln\left(\|\big(u(f,g)(\cdot,T),\
u(f,g)\big)\|_{H^2(\Omega)\times
{H^2(\Gamma\times(T-\delta_1,T+\delta_1))}}\right)\right|^{-1}.
\end{align*}
\end{pf}

\textbf{Acknowledgement}\\
\\
We are very grateful to Professor Michael Klibanov and Professor
Masahiro Yamamoto for pointing out the mistakes in our paper and
providing the guidance to correct these mistakes. Moreover,
Professor Michael Klibanov also provide his seminal papers on
Bukhgeim-Klibanov method \cite{K(78),K(79),B(80)}. And papers
\cite{I(1),Y(42)} modify the idea of the Bukhgeim-Klibanov method
via applying a new Carleman estimate for the parabolic operator of
Fursikov and Imanuvilov \cite{F(81)}. However, using the standard
Carleman estimate for the parabolic operator as in
\cite{K(78),K(79),B(80)}, one can also obtain H{\"o}lder stability.

\bibliographystyle{elsart-num-sort}
\bibliography{Recover}

\end{document}